\documentclass[12pt,dvips]{article}

\usepackage{amsmath,amsfonts,amssymb,amsthm,graphicx,verbatim,subfig}
\usepackage[all]{xy}

\ifx\pdftexversion\undefined

\usepackage[a4paper,colorlinks,link
=black,filecolor=black,citecolor=black,urlcolor=black,pdfstartview=FitH]{hyperref}
\else

\usepackage[a4paper,colorlinks,linkcolor=black,filecolor=black,citecolor=black,urlcolor=black,pdfstartview=FitH]{hyperref}
\fi

\font\sixbb=msbm6
\font\eightbb=msbm8
\font\twelvebb=msbm10 scaled 1095
\newfam\bbfam
\textfont\bbfam=\twelvebb \scriptfont\bbfam=\eightbb
                           \scriptscriptfont\bbfam=\sixbb
\def\bb{\fam\bbfam\twelvebb}
\newcommand{\Rea}{{\bb R}}
\newcommand{\Nat}{{\bb N}}
\newcommand{\Com}{{\bb C}}
\newcommand{\Int}{{\bb Z}}
\newcommand{\Rat}{{\bb Q}}
\newcommand{\FF}{{\bb F}}
\newcommand{\sgorf}{\overline{\FF}}
\newcommand{\KK}{{\bb K}}
\newcommand{\PP}{{\bb P}}

\newtheorem{theorem}{\bf Theorem}[section]
\newtheorem{claim}[theorem]{\bf Claim}

\newtheorem{proposition}[theorem]{\bf Proposition}
\newtheorem{corollary}[theorem]{\bf Corollary}
\newtheorem{observation}[theorem]{\bf Observation}
\newcommand{\enp}{\begin{flushright} $\Box$ \end{flushright}}
\newcommand{\beq}[0]{\begin{equation}}
\newcommand{\enq}[0]{\end{equation}}
\newcommand{\dn}{\Delta_{n-1}}

\newcommand{\cf}{{\cal F}}
\newcommand{\cl}{{\cal L}}

\newcommand{\rk}{\text{rk}\ }

\newcommand{\dnk}{\dn^{(k-1)}}

\newcommand{\cb}{{\cal B}}
\newcommand{\phat}{\widehat{\phi}}
\newcommand{\znk}{\Int_n^k}

\newcommand{\sg}{{\rm sgn}}

\newcommand{\namedref}[2]{\hyperref[#2]{#1~\ref*{#2}}}

\newcommand{\figureref}[1]{\namedref{Figure}{#1}}

\title{Sum Complexes - a New Family of Hypertrees}
\begin{document}
\author{N. Linial\thanks{Department of Computer Science, Hebrew University, Jerusalem 91904,
    Israel. e-mail: nati@cs.huji.ac.il~. Supported by ISF and BSF grants.} \and R. Meshulam\thanks{Department of Mathematics,
Technion, Haifa 32000, Israel. e-mail:
meshulam@math.technion.ac.il~. Supported by ISF and BSF grants. }
\and M. Rosenthal}

\maketitle
\pagestyle{plain}
\begin{abstract}
A $k$-dimensional hypertree $X$ is a $k$-dimensional complex on
$n$ vertices with a full $(k-1)$-dimensional skeleton and
$\binom{n-1}{k}$ facets such that $H_{k}(X;\Rat)=0$.  Here we
introduce the following family of simplicial complexes. Let $n,k$
be integers with $k+1$ and $n$ relatively prime, and let $A$ be a
$(k+1)$-element subset of the cyclic group $\Int_n$. The {\it sum
complex} $X_A$ is the pure $k$-dimensional complex on the vertex
set $\Int_n$ whose facets are $\sigma \subset \Int_n$ such that
$|\sigma|=k+1$ and $\sum_{x \in \sigma}x \in A$. It is shown that if
$n$ is prime then the complex $X_A$ is a $k$-hypertree for every
choice of $A$. On the other hand, for $n$ prime $X_A$ is
$k$-collapsible iff $A$ is an arithmetic progression in $\Int_n$.
\end{abstract}

\section{Introduction}
\ \ \ \ What is the high-dimensional analogue of a tree? Several
approaches to this question can be found in the literature. Here
we follow the lead of Kalai~\cite{kalai}. We start with some
standard notations. All simplicial complexes we consider $X$ have
$n$ vertices, and we always identify the vertex set of $X$ with
the cyclic group $\Int_n$. The number of $i$-dimensional faces of
$X$ is denoted by $f_i(X)$. We denote by $\dn$ the $(n-1)$-simplex
on the vertex set $\Int_n$
and by $\dn^{(i)}$ the $i$-dimensional skeleton of $\dn$.
A {\it $k$-hypertree} is a simplicial complex
$\dn^{(k-1)} \subset X \subset \dn^{(k)}$ such that $f_k(X)=\binom{n-1}{k}$ and
with a vanishing $k$-th rational homology $H_k(X;\Rat)=0$.
Throughout the paper we assume that $k+1$ is coprime to $n$.
For $a \in
\Int_n$, let $X_a$ be the following collection of subsets of
$\Int_n$:
$$X_a = \{\sigma \subset \Int_n:
|\sigma|=k+1~,~ \sum_{x\in \sigma} x =a\}~.$$ For a subset $A
\subset \Int_n$ of cardinality $k+1$, define the {\it Sum Complex}
$X_A$ by
$$X_A=\dnk \cup \left( \cup_{a \in A} X_a \right )~.$$
{\bf Example:} Let $n=7$, $k=2$ and $A=\{0,1,3\} \subset
\Int_7$. The $2$-dimensional complex $X_A$ (figure
\ref{fig:xa}) is obtained from the standard $6$-point
triangulation of the real projective plane $\Rea\PP^2$ on the
vertices $\{0,1,3,4,5,6\}$ (figure \ref{fig:rp2}) by
replacing the face $\{0,1,5\}$ with
the three faces $\{0,1,2\}$ , $\{0,2,5\}$ , $\{1,2,5\}$ , and
adding the faces $\{2,3,5\},\{0,2,6\}$ and $\{1,2,4\}$. $X_A$ is
clearly homotopy equivalent to $\Rea\PP^2$.

\begin{figure}

  \subfloat[The $6$-point triangulation of $\Rea\PP^2$]
  {\label{fig:rp2}
  \scalebox{0.7}{\input{twop1.pstex_t}}}
  \hspace{30pt}
  \subfloat[$X_A$ for $A=\{0,1,3\} \subset \Int_7$]
  {\label{fig:xa}
  \scalebox{0.7}{\input{twop2.pstex_t}}}
  \caption{}
  \label{figure1}
\end{figure}

In this paper we are concerned with topological and combinatorial
properties of $X_A$. Let $\FF$ be a field and let
$h_i(X_A;\FF)=\dim_{\FF} H_i(X_A;\FF)$. Since $X_A \supset \dnk$
it follows that $h_0(X_A;\FF)=1$ and $h_i(X_A;\FF)=0$ for $1 \leq
i \leq k-2$. Since $k+1$ is coprime to $n$, it follows that for
any $y \in \Int_n$, the number of $\sigma \subset \Int_n$ of
cardinality $k+1$ that satisfy $ \sum_{x \in \sigma} x=y$ is
$\frac{1}{n}\binom{n}{k+1}$. Therefore
$f_k(X_A)=\frac{k+1}{n}\binom{n}{k+1}=\binom{n-1}{k}$. The
Euler-Poincar\'{e} relation $\sum_{i \geq 0} (-1)^i
f_i(X_A)=\sum_{i \geq 0} (-1)^i h_i(X_A;\FF)$ then implies that
$h_{k-1}(X_A;\FF)=h_k(X_A;\FF)$. In the sequel we assume that the
characteristic of $\FF$ does not divide $n$.

Let $\omega$ be a fixed primitive $n$-th root of unity in the
algebraic closure $\sgorf$. For $x \in \Int_n$ let
$e(x)=\omega^x$. The $n \times n$ Fourier matrix $M$ over $\sgorf$
is given by $M(u,v)=e(-uv)$ for $u,v \in \Int_n$. For a subset $B
\subset \Int_n$ of cardinality $k+1$ let $M_{A,B}$ denote the
$(k+1)\times (k+1)$ submatrix of $M$ determined
 by $A$ and $B$. Let
$\cb_{n,k}$ denote the family of all $(k+1)$-element subsets of
$\Int_n$ that contain $0$.
\begin{theorem}
\label{hom}
\begin{equation}
\label{dimhom}
h_{k-1}(X_A;\FF)=h_k(X_A;\FF)=
\frac{1}{k+1}\sum_{B \in \cb_{n,k}} \dim\ker M_{A,B}~.
\end{equation}
\end{theorem}

The Fourier transform matrix $M=(M_{uv})$ of $\Int_n$ over
$\overline{\Rat}=\Com$ is given by $M_{uv}=\exp(-2 \pi i uv/n)$. A
classical result of Chebotar\"{e}v (see e.g. \cite{chebotarev})
asserts that if $n$ is prime then any square submatrix of $M$ is
nonsingular. Theorem \ref{hom} therefore implies
\begin{corollary}
\label{rational} If $n$ is prime then $X_A$ is a $k$-hypertree.
\end{corollary}

If $A$ is an arithmetic progression in $\Int_n$ then $M_{A,B}$ is
a Vandermonde matrix for all $B \in \cb_{n,k}$. Hence, by Theorem
\ref{hom}, $X_A$ is $\FF$-acyclic for any $\FF$ whose
characteristic is coprime to $n$. More is in fact true. Let
$\sigma$ be a face of dimension at most $k-1$ of a simplicial
complex $X$ which is contained in a {\it unique} maximal face
$\tau$ of $X$, and let $[\sigma,\tau]=\{\eta~:~ \sigma \subset
\eta \subset \tau \}$. The operation $X \rightarrow
Y=X-[\sigma,\tau]$ is called an {\it elementary $k$-collapse}. $X$
is {\it $k$-collapsible} if there exists a sequence of elementary
$k$-collapses $$X=X_1 \rightarrow X_2 \rightarrow \cdots
\rightarrow X_m=\{\emptyset\}~~.$$ Note that if  $\dn^{(k-1)}
\subset X \subset \dn^{(k)}$ is $k$-collapsible and
$f_k(X)=\binom{n-1}{k}$, then $X$ is $\Int$-acyclic.
\begin{theorem}
\label{collapse}
Let $n$ be a prime and
let $A$ be a subset of $\Int_n$ of cardinality $k+1$.
Then $X_A$ is $k$-collapsible iff $A$ is an arithmetic progression in $\Int_n$.
\end{theorem}

Theorems \ref{hom} and \ref{collapse} are proved in Sections
\ref{s:homology} and \ref{s:ap}. In Section \ref{s:examples} we
compute the homology of $X_A$ for $A=\{0,1,3\}$. We conclude in
Section  \ref{s:remarks} with some remarks concerning possible
extensions and open problems.

\section{Homology of $X_A$}
\label{s:homology}
\ \ \ \ We first recall some topological terminology (see e.g.
\cite{Munkres}). Let $X$ be a finite simplicial complex on the
vertex set $V$. For a set $S$ and a field $\KK$, let  $\cl(S,\KK)$
denote the $\KK$-linear space of all $\KK$-valued functions on $S$.
The space $C^m(X;\KK)$ of $\KK$-valued $m$-cochains of $X$ consists
of all functions $\phi \in \cl(V^{m+1},\KK)$ such that
$\phi(v_0,\ldots,v_m)={\rm sgn}(\pi)
\phi(v_{\pi(0)},\ldots,v_{\pi(m)})$ for any permutation $\pi$ on
$\{0,\ldots,m\}$, and such that $\phi(v_0,\ldots,v_m)=0$ if
$\{v_0,\ldots,v_m\}$ is not an $m$-dimensional simplex of $X$. (In
particular, $\phi(v_0,\ldots,v_m)=0$ if $v_i=v_j$ for some $i \neq
j$.) The coboundary operator $d_m:C^m(X;\KK) \rightarrow
C^{m+1}(X;\KK)$ is given by
$$d_m \phi (v_0,\ldots,v_{m+1})=
\sum_{i=0}^{m+1} (-1)^i \phi
(v_0,\ldots,\hat{v_i},\ldots,v_{m+1})~~.$$ Let $Z^{m}(X;\KK)=\ker
d_{m}$ denote the space of $m$-cocycles of $X$ over $\KK$ and
let $B^{m}(X;\KK)= {\rm Im ~} d_{m-1}$ denote the space of
$m$-coboundaries of $X$ over $\KK$. The $m$-dimensional cohomology
space of $X$ with coefficients in $\KK$ is
$$H^m(X;\KK)=\frac{Z^{m}(X;\KK)}{B^{m}(X;\KK)}~~.$$
Let $h^m(X,\KK)=\dim_{\KK}H^m(X;\KK)$. Then
$h^m(X,\KK)=h^m(X,\FF)=h_m(X;\FF)$ for any algebraic extension $\KK$
of $\FF$. In order to establish Theorem \ref{hom} we may therefore
assume that $\FF$ already contains a primitive $n$-th root of
unity $\omega$.

The Fourier transform of a function $\phi \in
\cl(\znk;\FF)$ is the function $\cf(\phi) =\phat \in \cl(\znk;\FF)$
given by
$$\phat(u_1,\ldots,u_k)=\sum_{(x_1,\ldots,x_k) \in \znk}
\phi(x_1,\ldots,x_k) e(-\sum_{j=1}^k u_j x_j)~~.$$ The Fourier
transform is an automorphism of $\cl(\znk;\FF)$.

The proof of Theorem \ref{hom} involves computing the image of
$H^{k-1}(X;\FF)$ under the Fourier transform. We first consider the
Fourier image of the $(k-1)$-coboundaries.
\begin{claim}
\label{bhat}
$$\cf(B^{k-1}(X_A;\FF))=
\{ g \in C^{k-1}(X_A;\FF) : \text{support}(g) \subset
\znk-(\Int_n-\{0\})^k \}~~.
$$
\end{claim}
\noindent {\bf Proof:} Let $\psi \in C^{k-2}(X_A;\FF)$. Then
$$\widehat{d_{k-2} \psi}(u_1,\ldots,u_k)=
\sum_{(x_1,\ldots,x_k) \in \znk} d_{k-2} \psi(x_1,\ldots,x_k)
e(-\sum_{j=1}^k u_j x_j)=$$ $$ \sum_{(x_1,\ldots,x_k) \in \znk}
\bigl( \sum_{i=1}^k (-1)^{i+1} \psi(x_1,\ldots,\hat{x_i},\ldots,x_k)
\bigr)
 e(-\sum_{j=1}^k u_j x_j)=
$$ $$
\sum_{i=1}^k (-1)^{i+1} \sum_{x_i} e(-u_i x_i) \sum_{x_1,\ldots
,\hat{x_i},\ldots,x_k} \psi(x_1,\ldots,\hat{x_i},\ldots,x_k)
e(-\sum_{j \neq i} u_j x_j) =$$
$$
n \sum_{i=1}^k (-1)^{i+1} \delta(0,u_i)
\sum_{x_1,\ldots
,\hat{x_i},\ldots,x_k} \psi(x_1,\ldots,\hat{x_i},\ldots,x_k)
e(-\sum_{j \neq i} u_j x_j)~~$$
where $\delta(0,u_i)=1$ if $u_i=0$ and is zero otherwise.
Therefore
$$\cf(B^{k-1}(X_A;\FF)) \subset \{ g \in C^{k-1}(X_A;\FF) : \text{support}(g)
\subset \znk-(\Int_n-\{0\})^k \}~~.$$ Equality follows since both spaces
have dimension $\binom{n-1}{k-1}$ over $\FF$. {\enp}

We next study the Fourier image of
the $(k-1)$-cocycles of $X_A$.
Fix a $\phi\in C^{k-1}(X_A;\FF)$. For $a
\in \Int_n$ define a function $f_a \in \cl(\znk;\FF)$ by
$$f_a(x_1,\ldots,x_k)=d_{k-1}\phi
\bigl(a-\sum_{i=1}^{k} x_i,x_1,\ldots,x_k\bigr)=
$$ $$\phi(x_1,\ldots,x_k)+ \sum_{i=1}^k (-1)^i
\phi\bigl(a-\sum_{j=1}^{k}
x_j,x_1,\ldots,\hat{x_i},\ldots,x_k\bigr)~~.$$ Let $T$ be the
automorphism of $\znk$ given by
$$
T(u_1,\ldots,u_k)=(u_2-u_1,\ldots,u_{k}-u_1,-u_1)~.$$
Then $T^{k+1}=I$ and for $1 \leq i \leq k$
$$
T^i(u_1,\ldots,u_k)=(u_{i+1}-u_i,\ldots,u_k-u_i,-u_i,u_1-u_i,
\ldots,u_{i-1}-u_i).$$
\begin{claim}
\label{ftfa} Let $u=(u_1,\ldots,u_k) \in \znk$. Then
\begin{equation}
\label{faone}
 \widehat{f_a}(u)=\phat(u)+ \sum_{i=1}^{k}(-1)^{ki}
e(-u_i a) \phat (T^{i} u)~~.
\end{equation}
\end{claim}
\noindent
{\bf Proof:} For $1 \leq i \leq k$ let
$\psi_i \in \cl(\znk,\FF)$ be given by
$$
\psi_i(x_1,\ldots,x_k)=
\phi\bigl(a-\sum_{j=1}^{k}
x_j,x_1,\ldots,\hat{x_i},\ldots,x_k\bigr)~~.$$
Then
$$
\widehat{\psi_i}(u)=
\sum_{(x_1,\ldots,x_k) \in \znk}
\phi\bigl(a-\sum_{j=1}^{k}
x_j,x_1,\ldots,\hat{x_i},\ldots,x_k\bigr)
e(-\sum_{j=1}^k u_j x_j)~~. $$
Substituting
$$
y_j=\left\{
\begin{array}{ll}
    a-\sum_{\ell=1}^k x_{\ell} & j=1 \\
    x_{j-1} & 2 \leq j \leq i \\
    x_j & i+1 \leq j \leq k
\end{array}
\right.
$$
it follows that
$$\sum_{j=1}^k u_j x_j = (a-y_1) u_i +\sum_{j=2}^i (u_{j-1}-u_i)y_j
+\sum_{j=i+1}^k (u_j-u_i)y_j~.$$
Therefore
$$
\widehat{\psi_i}(u)=e(-u_ia)
\sum_{y=(y_1,\ldots,y_k)\in \znk}
\phi(y)
e(u_iy_1-\sum_{j=2}^i(u_{j-1}-u_i)y_j-\sum_{j=i+1}^k (u_j-u_i)y_j)=$$
$$
e(-u_ia)\widehat{\phi}
(-u_i,u_1-u_i,\ldots,u_{i-1}-u_i,u_{i+1}-u_i,\ldots,u_k-u_i)=
$$
\begin{equation}
\label{psifour}
e(-u_ia) (-1)^{i(k-i)} \widehat{\phi}(T^i u)~.
\end{equation}
Now (\ref{faone}) follows from (\ref{psifour}) since
$f_a=\phi+\sum_{i=1}^k (-1)^i \psi_i$.
{\enp}
For $u \in \znk$ let $E_u=\{T^i u:0 \leq i \leq k\}$ and let
\begin{equation}
\label{phasez}
L_u=\bigcap_{a \in A}\{g \in \cl(E_u,\FF):g(u)+
\sum_{i=1}^{k}(-1)^{ki} e(-u_i a) g(T^{i} u)=0\}.
\end{equation}
Let $\phi \in Z^{k-1}(X_A;\FF)$. Then for all $a \in A$ and
$(x_1,\ldots,x_k) \in \znk$
$$f_a(x_1,\ldots,x_k)=d_{k-1}\phi
\bigl(a-\sum_{i=1}^{k} x_i,x_1,\ldots,x_k\bigr)=0~.$$
Eqn. (\ref{faone}) then implies that for all $a \in A$ and $u \in \znk$
$$\phat(u)+ \sum_{i=1}^{k}(-1)^{ki}
e(-u_i a) \phat (T^{i} u)=0~.$$
Writing $\phat_{| E_u}$ for the restriction of $\phat$ to $E_u$
we obtain
\begin{corollary}
\label{cor1}
Let $\phi \in C^{k-1}(X_A;\FF)$.
Then $\phi \in Z^{k-1}(X_A;\FF)$ iff
$\phat_{| E_u} \in L_u$ for all $u \in \znk$. {\enp}
\end{corollary}

Let the symmetric group $S_k$ act on $\znk$ by
$$\sigma\bigl( (u_1,\ldots,u_k)\bigr)=(u_{\sigma^{-1}(1)},\ldots,
u_{\sigma^{-1}(k)})~$$ and let $G_{n,k}$ denote the subgroup of
$\text{Aut}(\znk)$ generated by $T$ and $S_k$.
The subset $$D_{n,k}= \{(u_1,\ldots,u_k) \in (\Int_n-\{0\})^k:
u_i \neq u_j \text{~for~} i \neq j \}~~$$
is clearly invariant under $G_{n,k}$.
\begin{claim}
\label{gnk}
$~$
\\
(i) Let $\sigma \in S_k$ and $1 \leq i \leq k$.
Then $\eta=T^i  \sigma T^{-\sigma^{-1}(i)} \in S_k$ and
$\sg(\eta)=(-1)^{k(i+\sigma^{-1}(i))} \sg(\sigma)$.
\\
(ii) Any element of $G_{n,k}$ can be written
uniquely as $\sigma T^i$ where $\sigma \in S_k$
and $0 \leq i \leq k$. $G_{n,k}$ acts freely on $D_{n,k}$. 
\\
(iii) $L_u=L_{T^ju}$ for all $u \in D_{n,k}$ and $0 \leq j \leq k$.
\end{claim}
\noindent
{\bf Proof:} (i) For $1 \leq \ell \leq k$ let $\tau_{\ell} \in S_k$ be given by
$$
\tau_{\ell}(i)=\left\{
\begin{array}{ll}
    k-\ell+1+i & 1 \leq i \leq \ell-1 \\
    k-\ell+1 & i=\ell \\
    i-\ell & \ell+1 \leq i \leq k~.
\end{array}
\right.
$$
It can be checked that
$$\eta=T^i  \sigma T^{-\sigma^{-1}(i)}=\tau_{k-i+1}^{-1} \sigma \tau_{k-\sigma^{-1}(i)+1}~.$$
Noting that $\sg(\tau_{\ell})=(-1)^{k\ell+1}$ it thus follows that
$$\sg(\eta) =\sg(\sigma) \sg(\tau_{k-i+1}) \sg(\tau_{k-\sigma^{-1}(i)+1})= (-1)^{k(i+\sigma^{-1}(i))}\sg(\sigma)~~.$$
\\
(ii) It follows from (i) that
$$G_{n,k}=\{\sigma T^i: \sigma \in S_k~,~0 \leq i \leq k\}~.$$
Let $u=(u_1,\ldots,u_k) \in D_{n,k}$ and let
$v=(v_1,\ldots,v_k)= \sigma T^i u$. If $i \neq 0$ then
$$\sum_{j=1}^k v_j =\sum_{j=1}^k u_j  -(k+1) u_i \neq \sum_{j=1}^k u_j$$
and therefore $\sigma T^i u \neq u$. It follows that $G_{n,k}$
acts freely on $D_{n,k}$ and that the representation of an element of  $G_{n,k}$
as  $\sigma T^i$ is unique.
\\
(iii) Let $g \in L_u$ and $a \in A$. Then
$$
g(T^ju)+ \sum_{i=1}^k (-1)^{ik} e(-(T^j u)_i a)
g(T^{i+j}u)=$$
$$
g(T^ju)+\sum_{i=1}^{k-j} (-1)^{ik} e(-(u_{i+j}-u_j)a)g(T^{i+j}u)+ $$
$$
(-1)^{(k-j+1)k}e (u_ja)g(u)+ \sum_{i=k-j+2}^k (-1)^{ik}e(-(u_{i-k+j-1}-u_j)a)g(T^{i+j}u)=
$$
\begin{equation}
\label{gtu}
(-1)^{jk} e(u_j a) \bigr(g(u)+ \sum_{i=1}^k (-1)^{ik} e(-u_ia)g(T^iu)\bigl)=0.
\end{equation}
Hence $g \in L_{T^ju}$.
{\enp}
\noindent
{\bf Proof of Theorem \ref{hom}:}
Let $R \subset D_{n,k}$ be a fixed set
of representatives of the orbits of $G_{n,k}$ on $D_{n,k}$. Then
 $|R|=\frac{|D_{n,k}|}{|G_{n,k}|}=\frac{1}{k+1} \binom{n-1}{k}$.
Consider the mapping
$$\Theta: Z^{k-1}(X_A;\FF) \rightarrow \bigoplus_{u \in R} L_u$$ given by
$$\Theta (\phi)=\bigl( \widehat{\phi}_{|E_u}~:~u \in R \bigr)~~.$$
\begin{claim}
\label{injphi}
$$\ker \Theta =B^{k-1}(X_A;\FF)~~.$$
\end{claim}
\noindent
{\bf Proof:} $$\ker \Theta =
\{\phi \in  Z^{k-1}(X_A;\FF):\widehat{\phi}_{|E_u}=0
\text{~for~all~} u \in R \}=$$
$$\{\phi \in  Z^{k-1}(X_A;\FF):\widehat{\phi}(u)=0
\text{~for~all~} u \in D_{n,k}\}=B^{k-1}(X_A;\FF)$$
by Claim \ref{bhat}. {\enp}
\begin{claim}
\label{surjphi}
$\Theta$ is surjective.
\end{claim}
\noindent {\bf Proof:}
Let $(g_u : u \in R)
\in \bigoplus_{u \in R} L_u$. Define
$g \in C^{k-1}(X_A;\FF)$
by
$$
g(v)= \left\{
\begin{array}{ll}
    0 & v \not\in D_{n,k} \\
    \sg (\sigma) g_u(T^j u) & v=\sigma T^j u \text{~where~} u \in R.
\end{array}
\right.
$$
Clearly $\Theta(\cf^{-1}(g))= (g_u : u \in R)$.
To show that $\cf^{-1}(g) \in Z^{k-1}(X_A;\FF)$ it suffices
by Corollary \ref{cor1} to check that $g \in L_v$ for all $v \in \znk$.
If $v \not\in D_{n,k}$ then $g_{| E_v}=0$. Suppose then that
$v=\sigma T^j u \in D_{n,k}$ where $u \in R$ and $0 \leq j \leq k$. Combining Claim \ref{gnk}(i) and Eq. (\ref{gtu}) it follows that
$$g(v)+\sum_{i=1}^k (-1)^{ik} e(-v_i a) g(T^i v)=$$
$$ g(\sigma T^j u)+\sum_{i=1}^k (-1)^{ik}
e(-(\sigma T^j u)_i a) g(T^i\sigma T^j u)=$$
$$\sg(\sigma) g_u(T^j u)+ \sum_{i=1}^k (-1)^{ik}
e(-(T^ju)_{\sigma^{-1}(i)} a) (-1)^{k(i+\sigma^{-1}(i))} \sg(\sigma)
g_u(T^{\sigma^{-1}(i)+j}u)=$$
$$\sg(\sigma) \bigl(g_u(T^ju)+ \sum_{i=1}^k (-1)^{ik} e(-(T^j u)_i a)
g_u(T^{i+j}u)\bigr)=$$
$$=(-1)^{jk}\sg(\sigma)e(u_j a) \bigr(g_u(u)+ \sum_{i=1}^k (-1)^{ik} e(-u_ia)g_u(T^iu)\bigl)=0.
$$ {\enp}
Claims \ref{injphi} and \ref{surjphi} imply that
\begin{equation}
\label{mainiso}
H^{k-1}(X_A,\FF) \cong \bigoplus_{u \in R} L_u~.
\end{equation}
For $u=(u_1,\ldots,u_k) \in D_{n,k}$ let
$B_u=\{0,u_1,\ldots,u_k\}$. Then $\dim L_u= \dim \ker M_{A,B_u}~$.
Combining (\ref{mainiso}) with Claim \ref{gnk}(iii) it thus follows that
$$h^{k-1}(X_A;\FF)=
\sum_{u \in R}\dim L_u =$$
$$\frac{1}{k+1} \sum_{u \in R} \sum_{j=0}^k \dim L_{T^j u}=
\frac{1}{k+1}\sum_{B \in \cb_{n,k}} \dim\ker M_{A,B}~. $$
 {\enp}

\section{When is $X_A$ collapsible?}
\label{s:ap}

In this section we prove Theorem~\ref{collapse}, so that in this
section $n$ is prime. We find it convenient to maintain the
vertices in a face sorted according to the order induced from
$\Nat$, and also refer to subsets of $\FF_n$ as sorted vectors and
not only as sets.

\subsection{Equivalence}
Let $\phi: \FF_n \rightarrow \FF_n$ be the linear map $\phi(x) =
\alpha x + \beta$. It is clear that the image of $X_a$ under
$\phi$ is $X_t$ where $t=\alpha a + (k+1) \beta$. We say that the
complexes $X_{a_0, \ldots,a_k}$ and $X_{b_0,\ldots,b_k}$ are  {\em
equivalent} iff there exist a permutation $\pi$ on $\{b_0, \ldots,
b_k\}$ and $\alpha, \beta$ s.t. $\pi(b_i)=\alpha a_i + (k+1)
\beta$ for every $ 0 \leq i \leq k$. Equivalent complexes are
clearly isomorphic.

It is an easy observation that $a_0, \ldots, a_k$ is an arithmetic
progression iff $X_{a_0, \ldots, a_k}$ is equivalent to the
complex $X_{0,\ldots, k}$. We show that $X=X_{0,\ldots, k}$ is
collapsible whence $X_{a_0, \ldots, a_k}$ is collapsible for $a_0,
\ldots, a_k$ an arithmetic progression.

\subsection{Proof of sufficiency}

To show that $X$ is collapsible we introduce an order $\prec_R$ by
which we remove the $k$-faces from $X$. We need first some
preliminary definitions. With every $k$-face $u \in X$ we
associate a vector $h(u)$ of dimension $\lceil \frac{k}{2}
\rceil$. The $i$-th coordinate in $h$ counts how many integers in
the interval $[u_{i}, u_{k-i}]$ do not belong to $\{u_{i}, \ldots,
u_{k-i}\}$. Namely, the $i$-th coordinate of $h(u)$ is:
$$
h_i(u):= u_{k-i} - u_i -(k-2i)
$$
Clearly $h_i(u)$ is non-increasing in $i$. For every two $k$-faces
$u,v \in X$ we say that $u \prec_L v$  if $h(u)$ is
lexicographically smaller than $h(v)$. When $h(u) = h(v)$ we say
that $u \equiv_L v$. It should be clear that $h$ is invariant
under set reversal i.e. $x \rightarrow n-x$. It is also invariant
under shifts that ``do not overflow'' in the obvious sense, but we
will not be using this fact. If $u \not \equiv_L v$ for some $u,v
\in X$, we denote by $\delta_L(u,v)$ the first index for which
$h(u)$ and $h(v)$ differ. Thus if $u \prec_L v$ and $\delta_L(u,v)
= i$ then $h_j(u) = h_j(v)$ for all $j < i$ and $h_i(u) < h_i(v)$.\\
For $i, j \in \FF_n$ it is convenient to define $\rho (i,j)$ as
$i-j$ if $i > j$ and as $j-i$ otherwise. This is extended as usual
to: $\rho (i,A) = \min \{ \rho (i,a) \mid a\in A \}$ and $\rho
(A,B) = \min \{ \rho (a,b)\mid a\in A , b\in B \}$.\\
If $u \in X_i$ and $v \in X_j$ we say that $u \prec_I v$  if $i$
is closer than $j$ to $\{0,k\}$, i.e., if $\rho (i, \{0,k\}) <
\rho (j, \{0,k\})$. We say that $u \equiv_I v$ when $\rho (i,
\{0,k\}) = \rho (j, \{0,k\})$, namely, $i=j$ or $i = k-j$ .
Letting  $i' = \rho (i, \{0,k\})$, it is clear that $u \prec_I v$
iff $i' < j'$. If $u \not \equiv_I v$, we denote by $\delta_I(u,v)
= \min\{i', j'\} = \rho (\{i,j\}, \{0,k\})$ .

\begin{figure}[htbp]
\begin{center}
\caption{The order of collapse determined by $\prec_R$}
\label{fig:order}
\input{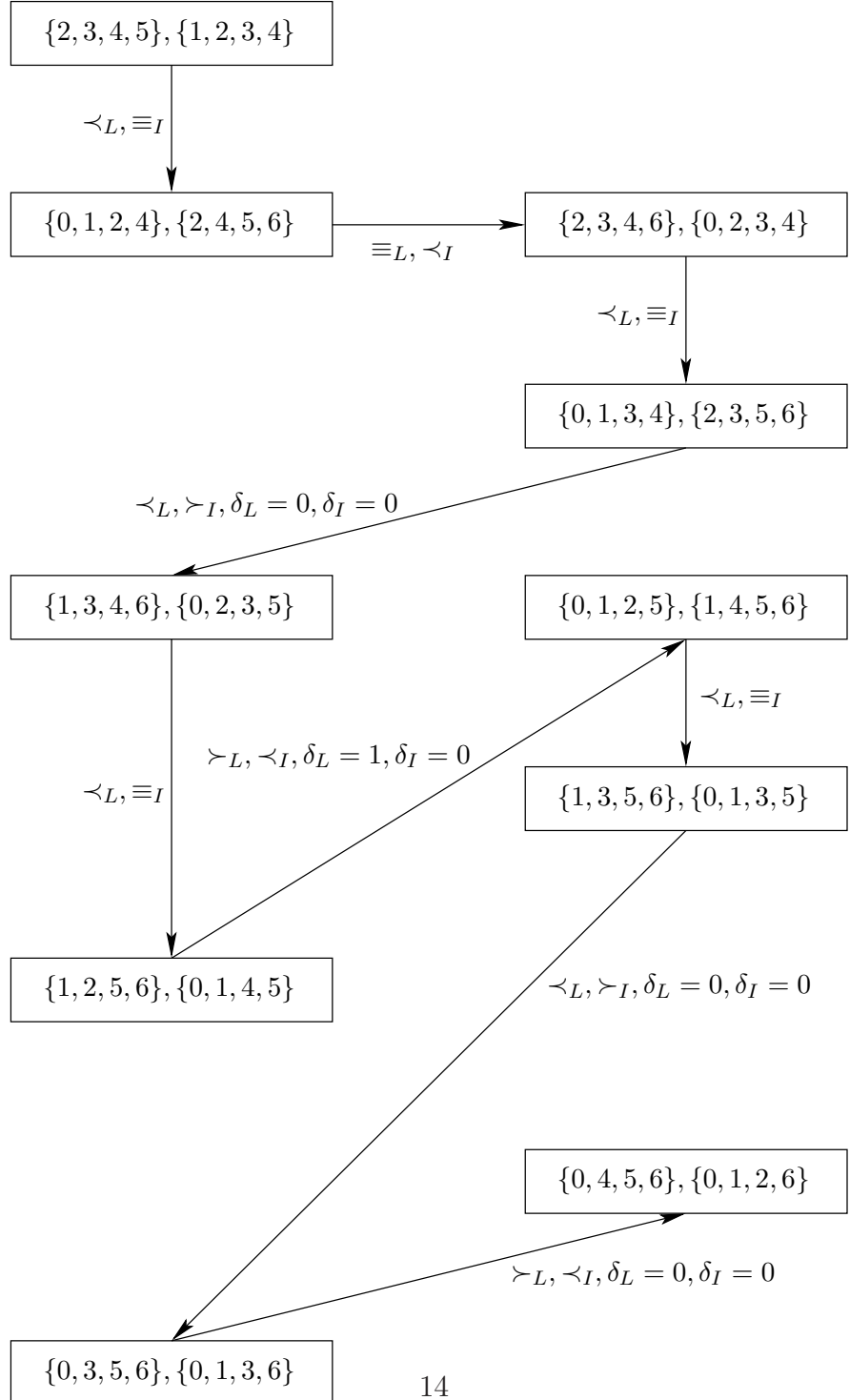}
\label{figure:order}
\end{center}
\end{figure}

We are now ready to define the relation $\prec_R$. This is done in
terms of the relations $\prec_L$ and $\prec_I$. To begin, $u
\equiv_R v$ iff $u \equiv_L v$ and $u \equiv_I v$. If $u \preceq_L
v$ and $u \preceq_I v$ and at least one inequality is proper, then
$u \prec_R v$. Finally, when $u \prec_L v$ and $u \succ_I v$, the
order $\prec_R$ is determined according to the smaller of
$\delta_I(u,v),\delta_L(u,v)$. Namely, if
$\delta_I(u,v)<\delta_L(u,v)$ then $u\succ_R v$ and if
$\delta_I(u,v) \geq \delta_L(u,v)$ then $u \prec_R v$.\\
To sum up, for $u,v \in X$:
\begin{enumerate}
\item If $u \equiv_I v$ and $u \equiv_L v$ then $u \equiv_R v$.
\item If $u \equiv_I v$ and $u \prec_L v$ then $u \prec_R v$.
\item If $u \prec_I v$, then $u \prec_R v$ unless
\begin{enumerate}
\item $u \succ_L v$ and
\item $\delta_L(u, v) \leq \delta_I(u,v)$
\end{enumerate}
In which case $u \succ_R v$.
\end{enumerate}
To clarify this definitions a little bit more, we present an
example from the complex $X_{0,1,2,3}$ over $\FF_7$. Let
$u=\{0,1,2,5\}$, $v=\{1,2,5,6\}$. The set $u$ has two missing
integers between $0$ and $5$ and no missing integers between $1$
and $2$, hence $h(u)=h(\{0,1,2,5\})=(2,0)$. Similarly
$h(v)=h(\{1,2,5,6\})=(2,2)$. Also, $u\prec_L v$ because $(2,0)$ is
lexicographically smaller than $(2,2)$. Furthermore,
$\delta_L(u,v)=1$ because the first coordinate the vectors differ
is the second coordinate (and we start indexing coordinates from
zero). Now $u\in {X_1}$ since $0+1+2+5\equiv 1 \mod 7$. Similarly
$v\in {X_0}$. We next calculate that $1' = 1 = \rho(1, \{0,7\})$
and $0' = 0$. Hence $v\prec_I u$ because $0' < 1'$, and
$\delta_I(u,v) = \min\{0,1\} = 0$. To recap, $u\prec_Lv$ and
$v\prec_I u$, so we turn to compare
$0=\delta_I(u,v)<\delta_L(u,v)=1$, it follows that in this case
the order $R$ is determined by $I$, hence $\{0,1,2,5\}=u \succ_R
v=\{1,2,5,6\}$. A full description of the order $R$ on
$X_{0,1,2,3}$ over $\FF_7$ is shown in \figureref{fig:param} and
\figureref{fig:order}:
\begin{figure}[htbp]
\caption{$X_{0,1,2,3}$ parameters over $\FF_7$}
\centering 
\begin{tabular}{ r  l  c  c} 
\\[1ex]
$\delta_L(x,y)$ & $h(x)$ & $i' = 0$ & $i' = 1$
\\
\hline \hline
&$(0,0)$ & $\{2,3,4,5\}, \{1,2,3,4\}$ & \\[-0.5ex]
\raisebox{1.5ex}{0} &$(1,0)$ & $\{0,1,2,4\}, \{2,4,5,6\}$ &
$\{2,3,4,6\}, \{0,2,3,4\}$ \\[-0.5ex]
\raisebox{1.5ex}{1}&$(1,1)$ & & $\{0,1,3,4\}, \{2,3,5,6\}$  \\[-0.5ex]
\raisebox{1.5ex}{0} &$(2,0)$ & $\{1,3,4,6\}, \{0,2,3,5\}$&
$\{0,1,2,5\}, \{1,4,5,6\}$  \\[-0.5ex]
\raisebox{1.5ex}{1}&$(2,1)$ & & $\{1,3,5,6\}, \{0,1,3,5\}$  \\[-0.5ex]
\raisebox{1.5ex}{1} &$(2,2)$ & $\{1,2,5,6\}, \{0,1,4,5\}$&   \\[-0.5ex]
\raisebox{1.5ex}{0} &$(3,0)$ & & $\{0,4,5,6\}, \{0,1,2,6\}$  \\[-0.5ex]
\raisebox{1.5ex}{1} &$(3,1)$ & $\{0,3,5,6\}, \{0,1,3,6\}$ & \\
\hline \\[4ex]
\end{tabular}\label{fig:param}
\end{figure}

A few words are in order about Figure~\ref{fig:param}. The rows
are sorted by the lexicographic order of $h(\cdot)$. The columns
on the right include all facets of $X$ sorted by value of $i'$.
Note that for each value of $h$ and each $i'$ there are two facets
that attain this pair of values. The leftmost column gives the
value of $\delta_L(x,y)$ for every two consecutive lines in the
table.

We now turn to show that $X$ can indeed be collapsed in the order
$\prec_R$. That is, for every $x \in X$ it is possible to apply an
elementary collapse step to $x$ if all the $\prec_R$-predecessors
of $x$ have already been collapsed. In order to show this, we need
to point out an free$(k-1)$-face that is contained in $x$. What we
will show is that for $x \in X_a$, the face $\hat x := x \setminus
\{x_a\}$ is free. (Note that since $x \in X = X_{0,\ldots,k}$,
there indeed must exist some $a \in \{0,\ldots,k\}$ s.t. $x \in
X_a$). It may be helpful to mention that $a$ plays a double role
here. It is an index in the vector $x$ as well as the sum of the
elements of $x$. Being free means that all the $k$-faces
containing $\hat x$, precede $x$ in the order $\prec_R$. A
$k$-face that contains $\hat x$ has the form $y^{(b)} := \hat x
\cup \{x_a+ (b-a)\}$ with $0\leq b \leq k$ and $b \neq a$.
Clearly, $y^{(b)}$ is a $k$-face in $X$ iff $x_a+ (b-a) \not\in
\hat x$. Also, in this case $y^{(b)} \in X_b$, as we assume below.

The proof that $y^{(b)}$ $\prec_R x$ has two cases:
\begin{enumerate}
\item We first consider the case where $y^{(b)} \succeq_I x$. Since
 $y^{(b)} \in X_b$ and $x \in X_a$, the meaning of $y^{(b)} \succeq_I x$
is that $b' \geq a'$. Therefore $\delta_I(y^{(b)}, x)$ which is
the smaller of $a'$ and $b'$ equals $a'$. This means that $b$ lies
between $a$ and $k-a$ (whether $a$ or $k-a$ is bigger is
immaterial here).
\begin{itemize}
\item Consequently, $x_a + (b-a)$ is in the interval $[x_{a'}, x_{k-a'}]$.
It follows that the first and last $a' -1$ elements of $x$ and
$y^{(b)}$ are identical. In particular, $h_i(y^{(b)}) = h_i(x)$
for $i < a'$.
\item  We recall that $y^{(b)}$ is created by removing $x_a$
from $x$ and replacing it by the term $x_a + (b-a)$. Thus the
interval $[y^{(b)}_{a'}, y^{(b)}_{k-a'}]$ is shorter than
$[x_{a'}, x_{k-a'}]$. It follows that the first coordinate where
$h(y^{(b)})$ and $h(x)$ differ is the $a'$-th coordinate, where
$h_{a'}(y^{(b)}) < h_{a'}(x)$. Consequently, $y^{(b)} \prec_L x$
and $a' = \delta_L(y^{(b)}, x)$.
\item If $y^{(b)} \equiv_I x$ then we are done, because we already know
that $y^{(b)} \prec_L x$. By definition of $\prec_R$ this yields
the desired conclusion $y^{(b)} \prec_R x$.
\item If $y^{(b)} \succ_I x$ then from the previous points we conclude that
$a' = \delta_L(y^{(b)}, x)= \delta_I(y^{(b)}, x)$. To sum up,
$y^{(b)} \prec_L x$ and $\delta_L(y^{(b)}, x)= \delta_I(y^{(b)},
x)$, which yields by definition, $y^{(b)} \prec_R x$, as claimed.
\end{itemize}
\item Now consider the case $y^{(b)} \prec_I x$. This means that $b' <
a'$. Therefore $b' = \delta_I(y^{(b)}, x)$. Consequently $b$ does
not lie between $a$ and $k-a$.
\begin{itemize}
\item It follows that $x_a + (b-a) \in [x_{b'}, x_{k-b'}]$.
Consequently, the first and last $b'+1$ elements of $x$ and
$y^{(b)}$ are identical. In particular, $h_i(y^{(b)}) = h_i(x)$
for $i \leq b'$. Thus $\delta_L(y^{(b)}, x) > b'$.
\item If $y^{(b)} \preceq_L x$ then $y^{(b)} \prec_R x$ and we are done.
\item If $y^{(b)} \succ_L x$ then from the previous points we
conclude that $b' = \delta_I(y^{(b)}, x) < \delta_L(y^{(b)},x)$.
Hence $y^{(b)} \succ_I x$ and $\delta_I(y^{(b)}, x) <
\delta_L(y^{(b)},x)$ . Again, by definition, $y^{(b)} \prec_R x$,
as claimed.
\end{itemize}
\end{enumerate}
This completes the proof that $X_{0,\ldots,k}$ is collapsible and
hence that $X= X_{a_0,\ldots,a_k}$ is collapsible whenever
$a_0,\ldots,a_k$ is an arithmetic progression.

\subsection{Proof of necessity}

We now turn to show that if $a_0,\ldots,a_k$ is not arithmetic,
then $X_{a_0,\ldots,a_k}$ is not collapsible. In fact we show that
in this case exactly $k+1$
elementary collapse steps can be carried out.\\
For $X \subseteq \FF_n$ we denote as usual by $X+a$ the $a$-shift
of $X$, namely, the set $\{x+a | x \in X\}$. We start with the
following simple observation.
\begin{observation} \label{arith_prog_under_shift}
Let $n$ be a prime. A subset $X \subsetneq \FF_n$ is an arithmetic
progression iff there is an element $l$ for which $|(X+l)
\setminus X|=1$.
\end{observation}

When is the $(k-1)$-face $x_1, \ldots, x_k$ free a free face? This
is the case iff, for each $k \ge i \ge 1$ the element $x_i +
\sum_{j=1}^k x_j$ belongs to the set $\{a_0,\ldots,a_k\}$. If $x_i
+ \sum x_j = a_l$ it means that $x_1, \ldots, x_k$ cannot be
extended to a $k$-face in $X_{a_l}$. This translates into a linear
system of equations in $x_1, \ldots, x_k$ whose matrix has $2$'s
along the main diagonal and $1$'s elsewhere. Such a matrix is
nonsingular, so the solution is unique. Also, all the $k$ terms
$x_i + \sum x_j$ are distinct, so the only choice we have in
constructing this linear system is which of the $k+1$ elements in
$\{a_0,\ldots,a_k\}$ to omit. There are $k+1$ such choices which
yields $k+1$ distinct collapse steps that can be carried out.

We now explicitly describe the $k+1$ collapse steps that can be
carried out. Each of these collapsible faces has the form
$x^{(t)}:=\{a_0+l_t,\ldots, a_k+l_t\} \in X_{a_t}$ for some $l_0,
\ldots l_k$

The condition $x^{(t)} \in X_{a_t}$ determines $l_t$ via $l_t =
\frac{a_t - \sum_{i=0}^k a_i}{k+1}$. We claim that the face
$y:=x^{(t)} \setminus \{a_t+l_t\}$ is free. The sum of $y$'s
elements is $-l_t$, so that for every $i\neq t$ we need to add the
term $\{a_i+l_t\}$ to $y$ in order to attain the sum $a_i$. This
is,
however, impossible since $\{a_i+l_t\}$ is a member of $y$.\\

We turn to show that after these first $k+1$ collapse steps are
carried out, there remain no free $(k-1)$-faces in $X$. In order
for a $(k-1)$-face $y$ to be free following the above collapses,
$y$ has to be contained in exactly one of these $k+1$ collapsed
faces. Since $\{a_0,\ldots,a_k\}$ is not an arithmetic
progression, by Observation~\ref{arith_prog_under_shift}, the
intersection of any two of the $x^{(t)}$ contains at most $k-1$
elements. In particular there is no $(k-1)$-face that they both
contain. Thus we have to consider only $(k-1)$-faces $y$ which are
contained in one of the $x^{(t)}$ and exactly one more $k$-face.

It follows that $y$ must be of the form $x^{(t)} \setminus \{a_j +
l_t\}$ for some $j$ and $t$. The sum of $y$'s elements is $a_t -
a_j - l_t$. If $y$ is contained as well in a $k$-face $z \in
X_{a_i}$, then necessarily $z_i = z = y \cup \{a_i - a_t + a_j +
l_t\}$. We are assuming that $y$ becomes free with the collapse of
$x^{(t)}$, so there must be exactly one index $i$ for which $z_i$
is a legal $k$-face different from $x^{(t)}$. It follows that
$x^{(t)}$ and $x^{(t)} + (a_j - a_t)$ must have $k$ elements in
common. Again by Observation~\ref{arith_prog_under_shift} this
means that the elements in $x^{(t)}$ form an arithmetic
progression, a contradiction.
The proof of Theorem \ref{collapse} is now complete.

\section{Example: Homology of  $X_{\{0,1,3\}}$}
\label{s:examples}
\ \ \ \
For a prime $p$ and an integer $n$ indivisible by $p$, let
$U_{p,n}$ be the group of $n$-th roots of unity in
$\overline{\FF_p}$.
\begin{proposition}
\label{except} Let $k=2$, $A=\{0,1,3\}$. Let $p$ be a prime and
suppose $n$ is coprime to $3p$. Then
$$h_1(X_A;\FF_p)=\frac{1}{3}|\bigl\{ \{u,v\} \subset U_{p,n}-\{1\}: u \neq
v \text{~and~} 1+u+v=0 \bigr\}|~.$$
\end{proposition}
\noindent {\bf Proof:} Let $B=\{0,k,\ell\}$ with $0<k<\ell<n$ and
let $u=\omega^{-k}, v=\omega^{-\ell}$. Then
$$\det M_{A,B}= \det \left[
\begin{array}{ccc}
1 & 1 & 1 \\
1 & u & v \\
1& u^3 & v^3
\end{array}
\right] =
$$
$$
uv^3-vu^3+u^3-u-v^3+v=(u-1)(v-1)(v-u)(u+v+1)~~.$$
It follows that
$$
\rk M_{A,B}=
\left\{
\begin{array}{ll}
        2 & 1+u+v=0  \\
        3 & \text{otherwise}.
\end{array}
\right. $$ Thus the Proposition follows directly from Theorem
\ref{dimhom}. {\enp}
\begin{corollary}
\label{specialp} Let $k=2$, $A=\{0,1,3\}$. Let $p$ be a prime and
suppose $n=p^m-1$ is coprime to 3. Then
$$h_1(X_A;\FF_p)=\left\{
\begin{array}{ll}
     \frac{n-1}{6} & p=2  \\
     \frac{n-2}{6} & p=3 \\
     \frac{n-4}{6} & p >3.
\end{array}
\right. $$
\end{corollary}
\noindent {\bf Proof:} Clearly $\FF_{p^m}^*=U_{p,n}$. Therefore,
by Proposition \ref{except}
$$h_1(X_A;\FF_p)=\frac{1}{6}|\bigl\{ u \in \FF_{p^m}^*-\{1\}:  -(1+u) \not\in
\{0,1,u\} \bigr\}|.$$
The Corollary now follows since
$$\bigl\{ u \in \FF_{p^m}^*-\{1\}:  -(1+u) \not\in
\{0,1,u\} \bigr\}=
\left\{
\begin{array}{ll}
     \FF_{2^m}^*-\{1\}  & p=2  \\
     \FF_{3^m}^*-\{\pm 1\}  &  p=3 \\
     \FF_{p^m}^*-\{\pm 1,-2,-\frac{1}{2}\}  & p >3.
\end{array}
\right. $$
{\enp}

\section{Concluding Remarks}
\label{s:remarks}
\ \ \ \
Theorem \ref{hom} provides an explicit description of the homology
of the sum complex $X_A$ over fields of
characteristic coprime to $n$. In particular, it follows
via Chebotar\"{e}v's Theorem that
if $n$ is prime then
$X_A$ is $\Rat$-acyclic, i.e. $X_A$ is a $k$-hypertree.
When $A$ is an arithmetic progression, $X_A$ was shown to be
$k$-collapsible, and in particular $\Int$-acyclic. One natural question
is whether there exist other $A$'s for which $X_A$ is $\Int$-acyclic.
Kalai's $k$-dimensional Cayley's formula \cite{kalai} suggests that
most $k$-hypertrees are not $\Int$-acyclic. Likewise we conjecture that $X_A$ is not
$\Int$-acyclic for most $(k+1)$-subsets $A \subset \Int_n$. One possible approach to the
question of $\FF_p$-acyclicity of $X_A$ for primes $p\nmid n$ is via the following reduction.
Let $S_{\FF}(A)$ be the $\FF$-linear space of polynomials in $\FF[x]$ spanned by the monomials
$\{x^a:a \in A\}$. Theorem \ref{hom} then implies that $X_A$ is $\FF_p$-acyclic iff
$\deg {\rm gcd}(f(x),x^n-1) \leq k$ for all $0 \neq f(x) \in S_{\overline{\FF_p}}(A)$.


\begin{thebibliography}{99}

\bibitem{kalai}
G. Kalai,
Enumeration of $\Rat$-acyclic simplicial complexes,
{\it Israel J. Math.} {\bf 45}(1983) 337--351.

\bibitem{Munkres}
J. Munkres, {\it Elements of Algebraic Topology}, Addison-Wesley,
1984.

\bibitem{chebotarev}
P. Stevenhagen and H. W. Lenstra, Chebotar\"{e}v and his density
theorem, {\it Math. Intelligencer} {\bf 18}(1996) 26--37.

\end{thebibliography}
\end{document}